\documentclass{article}
\usepackage[english]{babel}
\usepackage{graphicx} % Required for inserting images
\usepackage{multirow}%
\usepackage{amsmath,amssymb,amsfonts}%
\usepackage{amsthm}%
\usepackage{mathrsfs}%
\usepackage[title]{appendix}%
\usepackage{xcolor}%
\usepackage{textcomp}%
\usepackage{manyfoot}%
\usepackage{booktabs}%
\usepackage{algorithm}%
\usepackage{algorithmicx}%
\usepackage{algpseudocode}%
\usepackage{listings}%
\theoremstyle{definition}
\newtheorem{theorem}{Theorem}
\theoremstyle{definition}
\newtheorem{definition}{Definition}

\usepackage{nccmath}
\usepackage{mathtools}

\usepackage{caption}
 
\usepackage{subcaption}

\title{Optimal Marital Strategies: How Couples Develop Successful Interaction Styles}
\author{Micah Henson, Mark Kot, Ka-Kit Tung  
\\Department of Applied Mathematics, University of Washington }

\date{February 2024}

\begin{document}

\maketitle

\begin{abstract}
The study of marriage dynamics and of strategies to reduce the likelihood of divorce has been an important research area for many years. Gottman's \cite{gottman_roles_1993} research on successful marriages revealed three matched interaction styles: conflict-avoiding, validating, and volatile. There has, however, been little progress in explaining how couples develop these styles of interaction and why failure to do so leads to failed marriages. In this paper, we show that these interaction styles arise as solutions to an optimal control problem where the couples jointly maximize a common goal. The validating style arises when the benefit from achieving joint happiness is balanced by the emotional cost of adopting a particular style.  The ubiquitous conflict-avoider style arises naturally when the couple does not care about the cost. The volatile style is not an optimal solution, but volatile marriages may still be successful for couples with highly positive natural dispositions.  The problem of the spouses having different goals in marriage is relevant to marriage repair, and this problem will be studied in the next paper using differential game theory.
\end{abstract}

\section{Introduction}

A number of researchers have developed mathematical models to study marriage success and failure.  They have done so by using stochastic models \cite{barley_stochastic_2011, batabyal_dynamic_1999}, economic models \cite{cornelius_search_2003, huber_considering_1980, kippen_whats_2013, grossbard_economics_1993},  and dynamical systems \cite{matsumoto_love_2017, bae_chaotic_2015, satsangi_dynamics_2012, son_time_2011, ashraful_mahfuze_economics_2015, rey_jose-manuel_mathematical_2010}. 
 
 One of the first mathematical models of love was the Romeo-and-Juliet model of Strogatz \cite{strogatz_love_1988}. Strogatz intended his model to be an educational tool for students learning differential equations; it has, however, been a seminal work in this area. Rinaldi \cite{rinaldi_love_1998} extended Strogatz's model by including the initial appeal in a relationship, reciprocation of love, and forgetting. Rinaldi et al. later wrote a paper \cite{rinaldi_love_2010} and a book \cite{rinaldi_modeling_2016} detailing the dynamics of love and appeal. Sprott \cite{sprott_dynamical_2004}, in turn, expanded on the Romeo-and-Juliet model and showed that seemingly simple ordinary differential equations on love can have complex dynamics that reflect plausible real-life experiences. 

Becker \cite{becker_theory_1973} introduced concepts from economics to define a ``marriage market" that defines success in terms of the utility of each couple. Utility increases with the consumption of household goods, as well as with income, education, and physical beauty. Becker \cite{becker_theory_1974} introduced the component of marital ``caring" and quantified its effect in the ``marriage market." This work is highly influential in papers in economics studying love.

Our approach is instead based on research by Gottman and collaborators \cite{gottman_marital_1989, gottman-levenson, gottman_roles_1993}, who spent decades interviewing couples. In that research, he quantified spousal interactions during 15-minute interviews with the Rapid Couple Scoring System (RCISS). This system measures the number of positive comments minus negative comments during each turn at speech. From such data, Gottman identified three types of stable interaction-styles that predict successful marriages. These styles are called conflict-avoiding, validating, and volatile. A validating couple interacts in both the positive and negative range of emotion. A conflict-avoiding couple has little interaction in their negative range of emotions and only interacts when they feel positive. Finally, a volatile couple has passionate fights in their negative range of emotions but does not interact when they are positive. 

Cook et al. \cite{cook_mathematics_1995} conducted significant research that bridged the gap between psychology and mathematics. They wrote a set of difference equations for spousal behavior in an attempt to develop a mathematical theory behind the research of Gottman\cite{gottman_roles_1993}. In this work, the authors converted RCISS code into theoretical broken-stick functions representing each type of interaction style.

Tung \cite{tung_marriage_2007} rewrote the Cook et al. model as a coupled, continuous-time differential-equation model for the ``positivity'' of each spouse over time.  The marital influence of each spouse was modeled as a linear function, with a slope that differed depending on whether the spouse was positive or negative.  In this model, a positive or negative long-term (i.e., equilibrium) solution of the differential equation determined whether a marriage was successful or not.

In this paper, we turn Tung's model into an optimal-control problem in which a couple maximizes a shared goal to find the best strategy for their marriage. Our solution predicts the same linear functions with differing slopes that were previously used as influence functions by Cook et al. \cite{cook_mathematics_1995} and Tung \cite{tung_marriage_2007}.

Our work is novel in that very few researchers \cite{liu_optimal_2007, ashraful_mahfuze_economics_2015} have taken optimal-control approaches, and those that have focused on an approach rooted in economic concepts within the framework of a ``marriage stock." Instead, we base our optimal-control model on a psychological framework. Furthermore, this paper is different from previous approaches because it maximizes the happiness of a spouse based on social interactions within the marriage as opposed to optimizing marital capital based on utility \cite{liu_optimal_2007, ashraful_mahfuze_economics_2015}. Preliminary work applying optimal control to Tung's model was completed by Hemmi \cite{Hemmi2008} in an undergraduate thesis. 

In section \ref{model}, we introduce our model. In section \ref{method}, we present the essential definitions and theorems from optimal control theory that we use. We apply this theory to our model in section \ref{analysis}, and we then describe our numerical methods in section \ref{numericalMethod}. In section \ref{results}, we present our results, and we provide examples that help one interpret our results. Finally, in section \ref{discussion1}, we summarize our conclusions and describe some directions for future research.

\section{Model} \label{model}
Our main goal is to determine how couples evolve into their interaction styles. We propose a model where $x_1(t)$ and $x_2(t)$ represent the positivity of each spouse over time as measured, for example, by RCISS. These are our state variables. 

 The state equations,

\begin{subequations} 
 \begin{align}
 \label{statex}
     \frac{dx_1}{dt} = r_1(\bar{x}_1-x_1) + u_1x_2,\\ 
     \label{statey}\frac{dx_2}{dt} = r_2(\bar{x}_2-x_2) + u_2x_1,
 \end{align}
 \end{subequations}
account for the internal and external influences on each spouse. The state equations have the same form for each spouse, and so, without loss of generality, we now describe the equation for spouse 1.

The internal portion of the equation for spouse 1 consists of a natural disposition, $\bar{x}_1$, and a specific return rate, $r_1$, to that natural disposition. The external portion is a marital influence function,

 \begin{equation}
     I_1(x_2) = u_1x_2,
 \end{equation}
that describes the effect of spouse 2 on spouse 1. This influence function includes a control variable, $u_1$, that spouse 1
uses to adjust his or her response to spouse 2's positivity or negativity.

Spouse 1 chooses $u_1$ with the goal of maximizing the shared objective functional J, 
\begin{equation}
 \max_{\substack{0 \leq u_1 \leq u_{1,max} \\
  0 \leq u_2 \leq u_{2,max}}}[ \, J =
   \int_0^T F(t, x_1,x_2,u_1,u_2) \, dt \,].
\end{equation}
We assume $u_1\geq 0$ because a negative reaction to a negative spouse could unnaturally increase $dx_1/dt$. We impose $u_{1, max}$ because infinitely high interaction with one's spouse is unrealistic. Each spouse can only devote so much energy to responding to their spouse.

We assume that both spouses wish to maximize the same shared objective functional with integrand
\begin{equation} \label{obj}
    F= \alpha[x_1 -\frac{1}{2\epsilon}(u_1-u_0)^2] + (1-\alpha)[x_2 - \frac{1}{2\epsilon}(u_2-u_0)^2 ].
\end{equation}
The couple's goal could be to prioritize both spouses' happiness or, if one spouse needs more support, to prioritize that one spouse.
The parameter $\alpha$ accounts for the relative importance of each spouse in the function $F$.
If the goal is to prioritize spouse 1's happiness, $\alpha = 1$.
If the couple wishes to prioritize their joint or average happiness, $\alpha = 1/2$.
Finally, if the couple wants to prioritize spouse 2's happiness,  $\alpha = 0.$
The parameter $\alpha$ can have any value in the closed interval $[0,1]$.

The portion of equation (\ref{obj}) that begins with $1/(2\epsilon)$ is the emotional cost that spouse $i$ incurs while interacting with their spouse.  The parameter $u_0 \geq 0$ is the ideal, effortless, or natural level of interaction between the spouses. The parameter $\epsilon > 0$ then sets the emotional cost of interacting with (or of ignoring) one's spouse in a non-ideal way. As $\epsilon\to\infty$, there is no penalty for any interaction style, and the problem becomes linear. In contrast, for finite $\epsilon$, the problem is nonlinear. The optimal solution needs to balance happiness and the emotional cost required to achieve happiness, as measured by the objective functional.  

Other, more flexible payoffs do exist. One could, for example, proceed in analogy with Keyfitz \cite{Keyfitz1972TheMO} and consider a harmonic mean that automatically weighs the least happy spouse more heavily.

\section{Method: Optimal Control} \label{method}
We have formulated an optimal-control problem.
The study of optimal-control problems originated with Pontryagin's group \cite{pontry}.
We now briefly summarize some key concepts, definitions, and theorems from
optimal-control theory that we need.  See Lenhart and Workman \cite{lenhart},
Sethi \cite{sethi}, and other standard texts for more on these topics.

\begin{definition} We wish to solve a \textit{standard two-state-variable, two-control-variable optimal-control problem.}
That is, we wish to determine
    \begin{subequations}\label{OC problem}
\begin{align}
   \label{objfunDEF} \max_{u_1,u_2} \textbf{ }[\, J = \int_0^T F(t, x_1 (t), x_2 (t), u_1 (t), u_2 (t)) \, dt \,]
\end{align}
subject to the state equations.
\begin{align}
\dot{x}_1(t) &= f_1(t, x_1(t),x_2(t),u_1(t), u_2(t)),\\
\dot{x}_2(t) &= f_2(t, x_1(t),x_2(t),u_1(t), u_2(t)),
\end{align}
the control constraints
\begin{align}
u_1(t) &\in U_1 = [0, u_{1,max}],\\
u_2(t) &\in U_2 = [0, u_{2,max}],
\end{align}
and the initial conditions
\begin{align}
x_1(0) &= x_1^0,\; x_2(0) = x_2^0 .
\end{align}
\end{subequations}
We require that the functions $f_1$, $f_2$, and $F$ be continuously differentiable.
\end{definition}

\medskip

We next introduce an auxiliary function, the Hamiltonian, that plays a key role in optimal control theory.

\medskip

\begin{definition}
    
  The \textit{Hamiltonian} for the above problem is
  \begin{equation}
      H(t,x_1,x_2,u_1,u_2,\lambda_1, \lambda_2) = F(t,x_1,x_2,u_1,u_2) + \sum_{k=1}^2 \lambda_{k} f_k(t,x_1,x_2,u_1,u_2),
  \end{equation}
where the $\lambda_i(t)$, $i = 1,2$ are additional variables known as \textit{adjoint variables}.
\end{definition}
\medskip

The usual first step in solving optimal-control problem (\ref{OC problem}) is to satisfy a set of necessary conditions.  These necessary conditions are given by Pontryagin's maximum princple.

\begin{theorem}\label{PMP}
If $u_1^* (t)$, $u_2^* (t)$, $x_1^* (t)$, and $x_2^* (t)$ are optimal controls and responses for problem (\ref{OC problem}), there exist piecewise-differentiable adjoint variables, $\lambda_{1} (t)$ and $\lambda_{2} (t)$, such that
\begin{equation}\label{max condition}
    H(t, x_1^* , x_2^* , u_1 , u_2 , \lambda_{1}, \lambda_{2}) \leq H(t, x_1^* , x_2^* , u_1^* , u_2^* , \lambda_{1}, \lambda_{2})
\end{equation}
for all admissible controls, $u_1 \in U_1$ and $u_2 \in U_2$, and for all $t$, $0 \leq t \leq T$.
In addition, the adjoint variables must satisfy the adjoint equations
\begin{equation}
    \dot{\lambda}_{i} (t) = -\frac{\partial H}{\partial x_i} (t, x_1^* , x_2^* , u_1^* , u_2^*, \lambda_{1}, \lambda_{2} ),
\end{equation}
   and the terminal (transversality) conditions
   \begin{equation}
       \lambda_{i}(T) = 0 ,
   \end{equation}
for $i = 1,2.$
\end{theorem}

\medskip

It is useful at this point to note that our optimization problem now has two state variables,
two control variables, and two adjoint variables.  In addition to maximization
condition (\ref{max condition}) for the control variables, we have four differential
equation that describe the dynamics of the state variables and the adjoint variables. 
We also have initial conditions for the state variables and terminal conditions for the adjoint variables.

The most important part of the above theorem is, arguably, the existence of the adjoint variables.
For optimal controls and responses, there must exist adjoint variables that satisfy the above
conditions.  These adjoint variables help determine the optimal controls.

Another key part of the above theorem is the maximization of the Hamiltonian with respect
to the control variables.  If the optimal controls lie in the interior of their control
intervals, then
\begin{equation}
    \frac{\partial H}{\partial u_i} = 0,
\end{equation}
for $i = 1,2$, for the optimal controls.

Optimal controls may also occur, however, at the endpoints of their control intervals.  This is especially
common for Hamiltonians that are linear in the control variables, as when, for our problem, $\epsilon \rightarrow \infty$.
For such linear-control problems we choose the control variables
on the basis of switching functions.

\begin{definition}
For a Hamiltonian that is linear in both control variables, the \textit{switching functions} are
\begin{equation}
    \Psi_i(t) = \frac{\partial H}{\partial u_i},
\end{equation}
 for $i = 1,2$.  In other words, the switching functions are the coefficients of the control variables in the Hamiltonian.
\end{definition}

If the Hamiltonian is linear in its control variables, we choose each control variable
on the basis of the sign of its switching function,

    \begin{equation}
        u_i^* = \begin{cases}
            0, &\Psi_i < 0,\\
            u_{i,max},  &\Psi_i > 0,
        \end{cases}
    \end{equation}
for $i = 1,2$.  Typically, switches occur when $\Psi_i = 0$, which leads to \textit{bang-bang controls}. 
In addition, we may occasionally find cases where $\Psi_i = 0$
over an interval, leading to singular solutions, $u_i^* = u_{i,s}(t)$. 

The above conditions are usually enough to determine candidate solutions.
With these candidates in hand, we may then try to prove optimality using
a sufficiency theorem, such as that of Arrow and Kurz \cite{Arrow_Kurz_1977}.
This theorem is written in terms of the maximized or derived Hamiltonian.

\begin{definition}
The \textit{maximized Hamiltonian,} $M$, is
\begin{equation}
    M(t, x_1,x_2,\lambda_{1}, \lambda_{2}) = \max_{u_1 \in U_1, u_2 \in U_2} H(t, x_1,x_2, u_1, u_2, \lambda_{1}, \lambda_{2}).
\end{equation}
\end{definition}

\medskip

\begin{theorem}\label{pontyLD}
Let the control variables $u_1^* (t)$ and $u_2^* (t)$ and the corresponding responses and adjoint variables satisfy all of the
conditions of Theorem \ref{PMP} for $0 \leq t \leq T$.  Then $u_1^* (t)$ and $u_2^* (t)$ are optimal controls if
$M(t, x_1,x_2,\lambda_{1}, \lambda_{2})$ is concave with respect to $x_1$ and $x_2$ for all $t \in [0,T]$.
\end{theorem}

\medskip

The adjoint variables have a traditional economic interpretation.
If the objective functional  in equation (\ref{objfunDEF}) represents profit or cost, the adjoint variables represent the shadow prices that a player would be willing to pay for marginal increases in each of the state variables (Sethi \cite{sethi}). In our context, however, the adjoint variables represent emotional costs. The Hamiltonian also has a traditional economic interpretation. It is the sum of the rate of change of assets and dividends accruing from the state variable (Dorfman \cite{Dorfman_1969}). In our case, the assets are the happiness that a spouse or partner accrues from their marriage.

\section{Analysis} \label{analysis} 
There are some situations where our optimal control problem becomes linear and exact analytic solution can be found. In other situations, the nonlinear problem is solved using both perturbation and numerical methods.
\subsection{Case: $\epsilon \rightarrow \infty$}
First, we consider the case when $\epsilon \rightarrow \infty$. We now apply our necessary conditions to the problem of maximizing the objective functional J, 
\begin{subequations}\label{eq:problem}
\begin{align} 
    \max_{\substack{0 \leq u_1 \leq u_{1,max} \\
  0 \leq u_2 \leq u_{2,max}}}  [J &= \int_0^T \alpha x_i + (1-\alpha)x_j  dt],
\end{align}
subject to the state equations
\begin{align}
\dot{x_1} &= r_1(\bar{x}_1 - x_1) + u_1x_2, \\
\dot{x_2} &= r_2(\bar{x}_2 - x_2) + u_2x_1,
\end{align}
and the initial conditions
\begin{align}
x_1(0) &= x_1^0, x_2(0) = x_2^0.
\end{align}
\end{subequations}

The Hamiltonian is
\begin{equation}
    H=  \alpha x_1 + (1-\alpha)x_2  + \lambda_{1}[r_1(\bar{x}_1 - x_1) + u_1x_2] + \lambda_{2}[r_2(\bar{x}_2 - x_2) + u_2x_1], 
\end{equation}
and the switching functions are
\begin{equation}
    \frac{\partial H}{\partial u_i} = \Psi_i = \lambda_{i}x_j,
\end{equation}
 for $i,j = 1,2$, $i\neq j$.  The adjoint variables, in turn, satisfy the adjoint equations
\begin{subequations}\label{linadj}
\begin{align}
        \dot{\lambda}_{1} &= -\frac{\partial H}{\partial x_1} = -\alpha + \lambda_{1}r_1 - \lambda_{2}u_2^*,\\
    \dot{\lambda}_{2} &= -\frac{\partial H}{\partial x_2} = -(1-\alpha) - \lambda_{1}u_1^* + \lambda_{2}r_2.
\end{align}
\end{subequations}
and the terminal conditions
\begin{equation}
\lambda_1(T) = 0,  \lambda_2(T) = 0.
\end{equation}

The range of the parameters where the controls are singular (where $\Psi_i=0$  over an entire interval) is quite small. 
(See Appendix B for further details.)  We therefore focus our attention on piecewise-constant bang-bang controls that satisfy
\begin{equation} \label{phi}
    u_i^*= \begin{cases}
            0, \textbf{ } \lambda_i x_j < 0,\\
            u_{i,max}, \textbf{ }\lambda_i x_j > 0,
        \end{cases}
\end{equation}
for $i,j = 1,2$, $i \neq j$.

\begin{figure}[ht!]
         \centering
    \includegraphics[width=\textwidth]{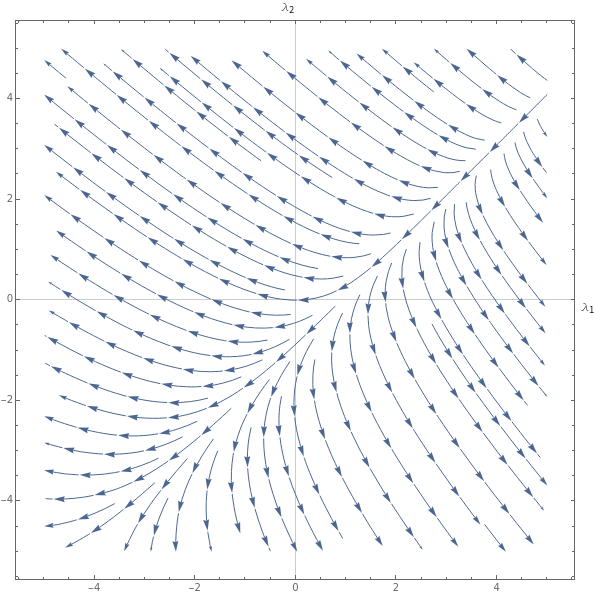}
         \caption{Sample phase-portrait for the adjoint space for parameters $r_1 = 0.6, r_2 = 0.6, \alpha = 1, \bar{x}_1 = 0.2, \bar{x}_2 = -0.4, u_1^* = 0.5,$ and $u_2^* = 0.5$.}
         \label{fig:adj_pp}
     \end{figure}

For each allowable combination of $u_1^*$ and $u_2^*$, there is only one orbit of adjoint system (\ref{linadj}) that enters the origin, and this orbit lies in the first quadrant of the $( \lambda_1, \lambda_2)$ phase-plane. (See Figure \ref{fig:adj_pp} for an example and Appendix \ref{proof} for
further details and analyses.)  Moreover, orbits outside the first quadrant are unable to enter
the first quadrant.   It quickly follows that $\lambda_1(t) > 0$ and $\lambda_2(t) > 0$ for all $0 \le t < T$. 

For each spouse, the optimal control here is thus the conflict-avoiding interaction style that Gottman observed in many successful marriages.  Each partner avoids conflict by ignoring their spouse when their spouse is negative. They do interact when their spouse is positive.  There is, of course, an emotional cost associated with ignoring one's spouse when he or she is upset, but each individual here ignores that cost. 

\subsection{Ubiquity of the conflict-avoider style}
A general case for which the conflict-avoider style arises as an optimal marital strategy is when the goals of the couple does not involve their happiness.  This could be the case where the couples, for example, are optimizing their financial well-being in marriage.  In this case, $F$ is independent of the $x's$. Another situation where the conflict avoider style appears as an optimal solution is when one of the spouses is not optimizing his own happiness--he could be the altruistic type. He is instead maximizing the happiness of the other spouse.  This could be the case when $\alpha=0$ or $\alpha=1$.

\subsection{The nonlinear problem}

For the complete model from section \ref{model}, the Hamiltonian is
\begin{equation}
    \begin{split}
        H =  &\alpha[x_1 -\frac{1}{2\epsilon}(u_1-u_0)^2] + (1-\alpha)[x_2 - \frac{1}{2\epsilon}(u_2-u_0)^2 ] \\
        &+ \lambda_{1}[r_1(\bar{x}_1 - x_1) + u_1x_2] + \lambda_{2}[r_2(\bar{x}_2 - x_2) + u_2x_1].
    \end{split}
\end{equation}
We now find the optimal controls, $u_1^*$ and $u_2^*$, by setting the derivatives of the Hamiltonian with respect to $u_1$ and $u_2$ equal to zero, 
\begin{subequations}\label{intCont}
\begin{align}
    \frac{\partial H}{\partial u_1} &= -\frac{\alpha}{\epsilon}(u_1-u_0) + \lambda_{1}x_2 = 0, 
    \\
    \frac{\partial H}{\partial u_2} &= -\frac{(1-\alpha)}{\epsilon}(u_2-u_0) + \lambda_{2}x_1 = 0. 
    \end{align}
\end{subequations}

If $\alpha = 0$, spouse 1's optimal control is as defined in equation (\ref{phi}). Likewise, if $\alpha = 1$,  spouse 2's optimal control is given by equation (\ref{phi}).  Otherwise,
\begin{subequations} \label{nphis}
\begin{gather}
    u_1^* =  \frac{\epsilon}{\alpha}\lambda_{1}x_2 + u_0,
    \\
    u_2^* =  \frac{\epsilon}{1-\alpha}\lambda_{2}x_1 + u_0.
\end{gather}
\end{subequations}

Solving the maximized system of equations, 
\begin{subequations}
    \begin{align}
    \dot{x}_1^* &= r_1(\bar{x}_1 - x_1^*) + u_1^* x_2^* \\
    \dot{x}_2^* &= r_2(\bar{x}_2 - x_2^*) + u_2^* x_1^* \\
    \dot{\lambda}_{1} &= -\alpha + \lambda_{1}r_1 -\lambda_{2}u_2^* \label{adj1}\\
     \dot{\lambda}_{2} &= -(1-\alpha) -\lambda_{1} u_1^* +\lambda_{2}r_2, \label{adj2}
    \end{align}
\end{subequations}
with initial conditions  $x_1(0) = x_1^0, x_2(0) = x_2^0$ and terminal conditions $\lambda_{1}(T) = 0$, $\lambda_{2}(T) = 0$, now gives us the solution to our problem.

For the nonlinear problem, we use perturbation theory to find analytical solutions for our optimal control problem for small $\epsilon$. The finite $\epsilon$ case will be dealt with numerically later in section \ref{results}.  The method has a long and storied history in applied mathematics, dating back to Henry Poincare in 1886 \cite{Murray_1984}. The procedure presented here is adapted from \cite{Bender_Orszag_2010}.

Consider the differential equation
\begin{equation}\label{genDE}
    \frac{dx}{dt} + p(t)x(t) = g(t)
\end{equation}
with initial condition $x(0) = A$.

The variable $x(t)$ can be expanded via a Taylor series in the parameter $\epsilon \rightarrow 0$,

\begin{equation} \label{tay}
    x(t;\epsilon) = x_0(t) + \epsilon x_1(t) + \epsilon^2 x_2(t) + ... = \sum_{n=0}^{\infty} \epsilon^n x_n(t).
\end{equation}
If we plug this Taylor expansion into equation (\ref{genDE}) and compare coefficients of $\epsilon_n$, we end up with an infinite system of differential equations that may (or may not) be easier to solve:
\begin{eqnarray}
    x_0'(t) + p(t)x_0(t) &=& g(t), \text{ } x_0(t) = A,\\
    x_1'(t) + p(t)x_1(t) &=& 0,  \text{ } x_1(t) = 0\\
    &\vdots&
\end{eqnarray}
Truncating the Taylor expansion of $x$ at some finite $n$ gives an approximate analytical solution to differential equation (\ref{genDE}).

\subsubsection{Zeroth order}
We now look at terms with no $\epsilon$ parameter for this order. This gives us a system, 
\begin{subequations}\label{zerosys}
    \begin{align} 
         \dot{x}_{1,0} = r_1(\bar{x}_1 - x_{1,0}) + u_0x_{2,0}, \label{zerothA} \\
          \dot{x}_{2,0} = r_2(\bar{x}_2 - x_{2,0}) +  u_0x_{1,0}, \label{zerothB}\\
           \dot{\lambda}_{1,0} = -\alpha + \lambda_{1,0}r_1 - \lambda_{2,0}u_0, \\
         \dot{\lambda}_{2,0} =-(1-\alpha) - \lambda_{1,0}u_0 + \lambda_{2,0}r_2
    \end{align}
\end{subequations}
for obtaining a solvable, analytical solution. 

The zeroth-order system of equations is now a system of three, uncoupled pairs of equations that can each be solved. For the zeroth-order state equations, we get
\begin{equation}
    \begin{bmatrix}
        x_{1,0} \\x_{2,0}
    \end{bmatrix} =
    a_1\begin{bmatrix}
        u_0\\r_1 + \eta_1
    \end{bmatrix}e^{\eta_1t}  + a_2\begin{bmatrix}
        u_0\\r_1 + \eta_2
    \end{bmatrix}e^{\eta_2t} + 
    \begin{bmatrix}
        \frac{\bar{x}_1r_1r_2 + \bar{x}_2r_2u_0}{r_1r_2-u_0^2} \\ \frac{r_1r_2\bar{x}_2 + r_1\bar{x}_1u_0}{r_1r_2-u_0^2}      
    \end{bmatrix}
\end{equation}
where 
\begin{equation}
    \eta_1,\eta_2 = \frac{1}{2}[-(r_1+r_2)\pm \sqrt{(r_1+r_2)^2 - 4(r_1r_2-u_0^2)}],
\end{equation}
and initial conditions 
\begin{equation}
     \begin{bmatrix}
        x_{1,0}(0) \\x_{2,0}(0)
    \end{bmatrix} = \begin{bmatrix}
        x_1^0 \\x_2^0
    \end{bmatrix}.
\end{equation}

For the two sets of adjoint equations, we have
\begin{equation}
    \begin{bmatrix}
        \lambda_{1,0} \\\lambda_{2,0}
    \end{bmatrix} =
    b_1\begin{bmatrix}
        u_0\\r_1 - \rho_1
    \end{bmatrix}e^{\rho_1t}  + b_2\begin{bmatrix}
        u_0\\r_1 - \rho_2
    \end{bmatrix}e^{\rho_2t} + 
    \begin{bmatrix}
         \frac{u_0-u_0\alpha + \alpha r_2}{r_1r_2-u_0^2} \\ \frac{u_0\alpha +(1-\alpha)r_1}{r_1r_2-u_0^2}      
    \end{bmatrix},
\end{equation}
where
\begin{equation}
    \rho_{1,2} = \frac{1}{2}[r_1+r_2\pm \sqrt{(r_1+r_2)^2 - 4(r_1r_2-u_0^2)}],
\end{equation}
with terminal conditions
\begin{equation}
    \begin{bmatrix}
        \lambda_{1,0}(T) \\ \lambda_{2,0}(T)
    \end{bmatrix} = \begin{bmatrix}
        0\\0
    \end{bmatrix}.
\end{equation}
 
\subsubsection{First Order}
This system of equations comes from the coefficients of terms with $\epsilon$ raised to one,
\begin{subequations} \label{firstOrder}
    \begin{align}
         \dot{x}_{1,1} = -r_1 x_{1,1} + \frac{1}{\alpha}\lambda_{1,0}x_{2,0}^2 + u_0x_{2,1} \label{firstA} \\
          \dot{x}_{2,1} = -r_2x_{2,1} + \frac{1}{(1-\alpha)}\mu_{2,0}x_{1,0}^2 + u_0x_{1,1},  \label{firstB}\\
           \dot{\lambda}_{1,1} = \lambda_{1,1}r_1 - \frac{1}{(1-\alpha)}\lambda_{2,0}\mu_{2,0}x_{1,0} - u_0\lambda_{2,1}, \\
         \dot{\lambda}_{2,1} = -\frac{1}{\alpha}\lambda_{1,0}^2x_{2,0} - u_0\lambda_{1,1} + \lambda_{2,1}r_2
    \end{align}
\end{subequations}
It only depends on zeroth order terms. All of these equations can be written

\begin{subequations}
    \begin{align}
        \dot{\vec{x}} + A \vec{x} = \vec{{f}}(t) \\
        \dot{\vec{\lambda}} - A \vec{\lambda} = \vec{{g}}(t)
    \end{align}
\end{subequations}
where
\begin{equation}
    A = \begin{bmatrix}
        r_1 & -u_0 \\ -u_0 & r_2,
    \end{bmatrix}
\end{equation}
with conditions $\dot{\vec{x}}(0) = 0$ and $\dot{\vec{\lambda}}(T) = 0$. Each equation can be solved with any two-dimensional, linear, differential-equation method. The exact solution of the first-order scheme contains many terms and parameters, so we omit it. Higher-order perturbation systems can systematically be determined and solved.

Overall, using the perturbation scheme, we recover a validating interaction-style by considering small values of $\epsilon$. Taking into account the penalty for ignoring one's spouse differentiates the validating style from the avoider style marriages. As $\epsilon \rightarrow 0$, there is a massive cost of adopting the optimal solution, and we see in equation (\ref{nphis}) that the first term of these equations also goes to zero, and the only interaction between the spouses comes from $u_0$. The zeroth-order perturbation scheme verifies this. In equations (\ref{zerothA}) and (\ref{zerothB}), we see a fixed, positive response of $u_0$ between spouses. This is equivalent to a validating interaction-style. 

 The higher orders of the perturbation scheme show how the different marriage goals affect the relationship because of the presence of the $\alpha_i$ values. They are deviations from the validating marriage. The goal of the marriage determines precisely how significant that deviation is. In equation (\ref{firstOrder}), a spouse with a balanced marriage goal gets twice the increase from their spouse's positivity compared to the selfish spouse. 

If $\alpha = 0$, spouse 1 never becomes validating. If $\alpha = 1$, spouse 2 never becomes validating. The optimal strategy remains conflict-avoiding regardless of the values of $\epsilon$. We call this type of spouse an altruist. The altruist decides to accept the higher emotional cost that comes with ignoring their spouse. We explore how marriage goals affect how couples interact with each other through specific numerical examples. 

\section{Numerical Method} \label{numericalMethod}
We use numerical methods to verify our analytical results. However, finding numerical solutions comes with challenges because the state equations have initial conditions, while the adjoint equations have terminal conditions. This means we have a two-point, boundary-value problem.  We have built our solvers using the Forward-Backward Sweep (FBS) algorithm. 

FBS specifically solves problems in optimal control when the Maximum Principle formulates solutions. The version presented here is based on Lenhart and Workman\cite{lenhart}, which is based on Wolfgang \cite{wolfgang}. McAsey et al. \cite{convOC} proved the convergence of this method for optimal control problems. The steps of the algorithm are:
\begin{enumerate}
    \item Discretize time and initialize step size. Guess an initial control.
 \item Use the initial condition to solve the state system forward in time.
 \item Use the terminal condition to solve the adjoint system backward in time. 
 \item Update the control using the state and adjoint solutions previously calculated.
 \item Repeat this process until the differences of each variable between iterations are within a specified tolerance. 
\end{enumerate}
 
Any appropriate numerical method for
solving ordinary differential equations that return values at the specified grid points can be used in steps 2 and 3. We chose to use Runge-Kutta 4. More information on this method can be found in Butcher \cite{Butcher_1967}.

\section{Results} \label{results}
 
We examine the implications of this model through a few numerical examples to verify our analytical result that $\epsilon \rightarrow \infty$ gives rise to conflict-avoiding marriages. We also show examples that consider how changing $\epsilon$ at finite numbers effects the model.

\subsection{Example 1: $\epsilon \rightarrow \infty$}
\begin{figure}[ht!]
    \centering
    \includegraphics[width=\textwidth, =\textwidth, trim= 1cm 2cm 1cm 2cm]{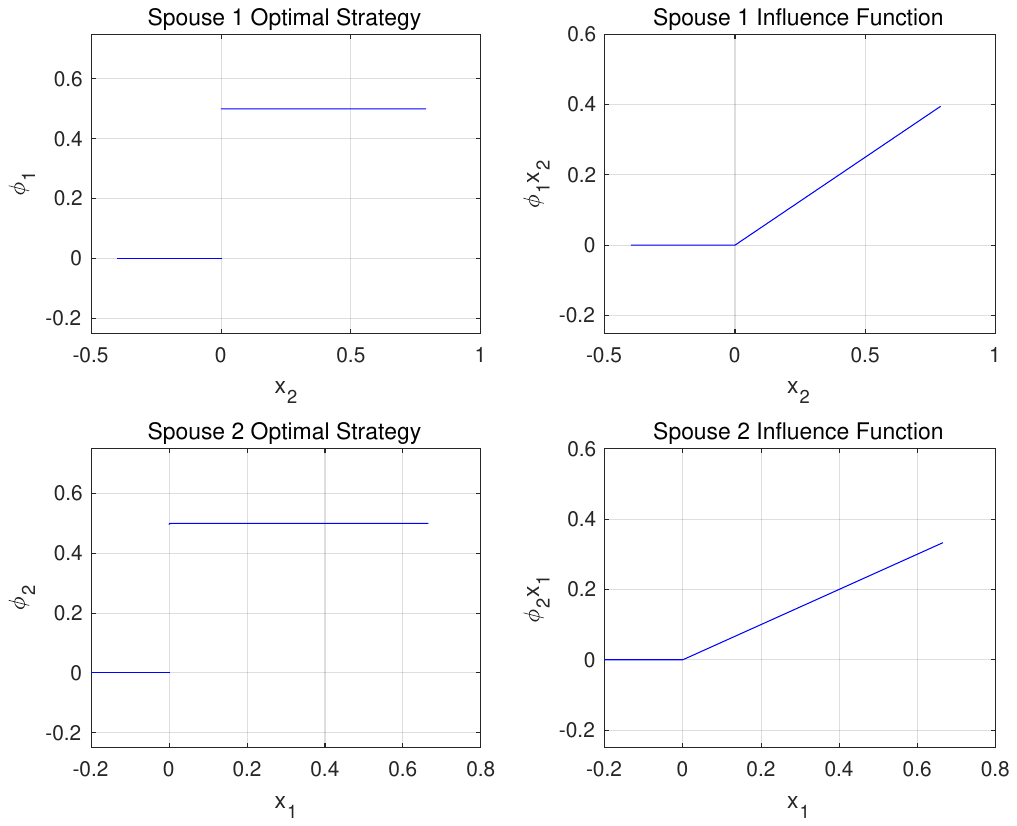}
    \caption{Sample optimal controls}
    \label{fig:ex1_oc}
\end{figure}

\begin{figure}[ht!]
    \centering
    \includegraphics[width=\textwidth]{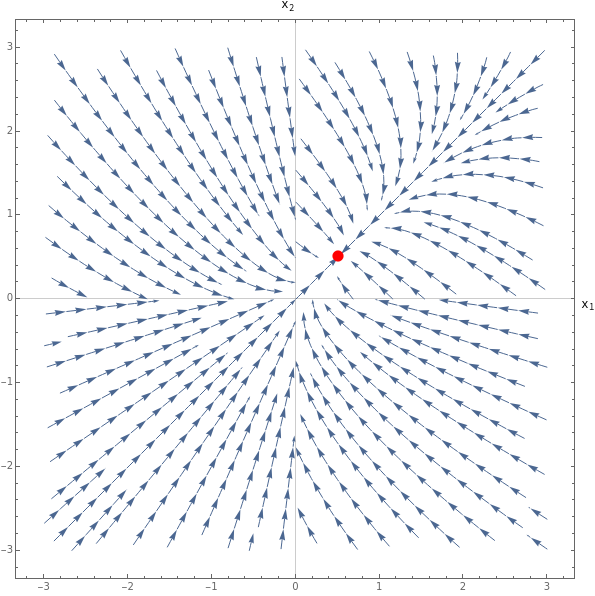}
    \caption{Example state-space with equilibrium point for a couple with positive natural disposition}
    \label{fig:ex1_state}
\end{figure}
We verify our analytical results with numerical results obtained by applying the method outlined in section \ref{numericalMethod}. Let us try an example with specific parameters. Let $u_{i,max} = 0.5$, $r_1 = 0.82 = r_2, \bar{x}_1 = 0.2, \bar{x}_2 = 0.4, x_i^0 = -\bar{x}_i, \alpha = 1$, and $T = 10$.

In figure \ref{fig:ex1_oc}, each spouse's control variable $u_i^*$ switches from $0$ to $1$ as their partner's positivity switches from negative to positive. These are conflict-avoiding interaction-styles. They are of exactly the same form as the theoretical influence functions in Tung.  In figure \ref{fig:ex1_state}, we have the phase portrait for the couple in this example. We see that the couple heading to an equilibrium point in the first quadrant. This means that the couple in this marriage will ultimately be successful in the sense that they are happier in marriage than if they are single, or stable, meaning that they are both happy. 

\subsection{Example 2: $\epsilon$ is finite}

\begin{figure}[ht!]
    \centering
    \includegraphics[width=\textwidth, =\textwidth, trim= 1cm 2cm 1cm 2cm]{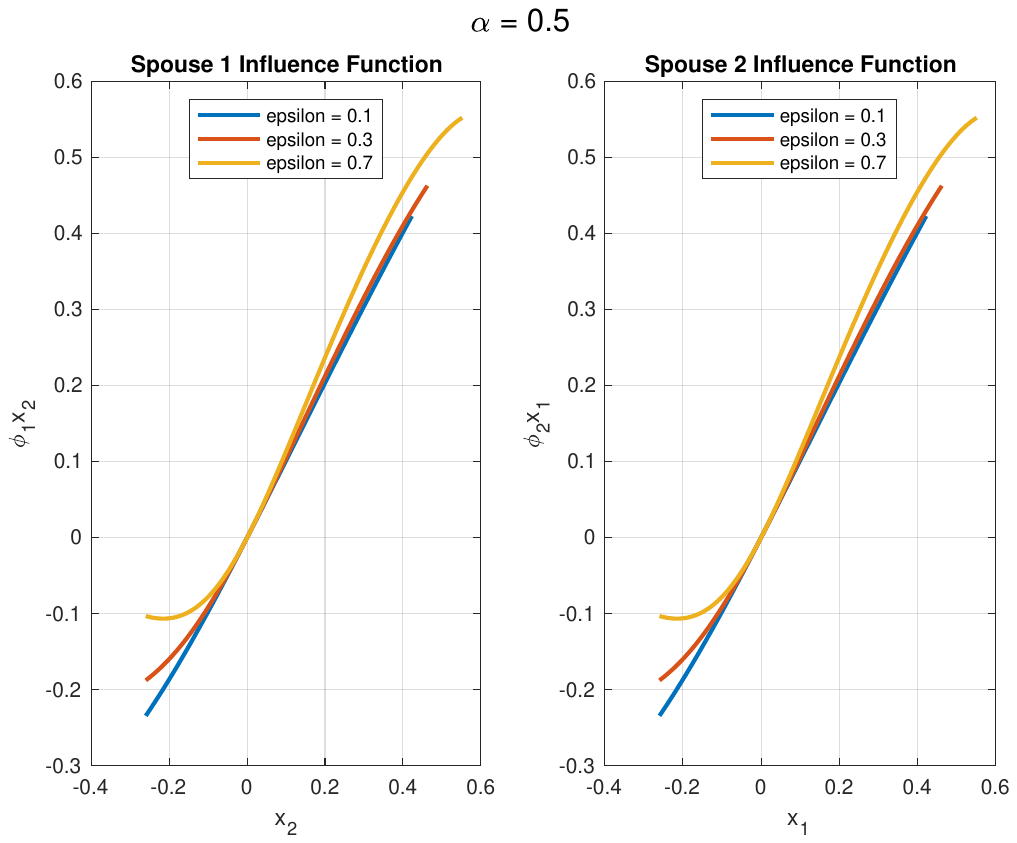}
    \caption{Sample Influence Functions for varying $\epsilon$ values and fixed $\alpha = 1/2$.}
    \label{fig:ex2}
\end{figure}

For parameters
$r_1 = 0.82 = r_2, \bar{x}_1 = 0.26 =\bar{x}_2, x_i^0 = -\bar{x}_i$, $\alpha = 1/2$ and $T = 3$, we plot the influence function for different values of $\epsilon$. The results can be found in figure \ref{fig:ex2}. This figure verifies what we saw in the analysis, specifically in the perturbation solution. We can see that both spouses have adapted a validating interaction style for $\alpha = 0.5$ and $\epsilon$ small. This figure also verifies that $\epsilon$ increasing results in deviation from a validating marriage style, as was determined in the first order perturbation scheme. As $\epsilon$ increases further, the influence function looks less like a validating style and moves closer to a conflict-avoiding style.

\subsection{Example 3: Changing $\alpha$ for a fixed and finite $\epsilon$}

\begin{figure}[ht!]
    \centering
    \includegraphics[width=\textwidth, =\textwidth, trim= 1cm 2cm 1cm 2cm]{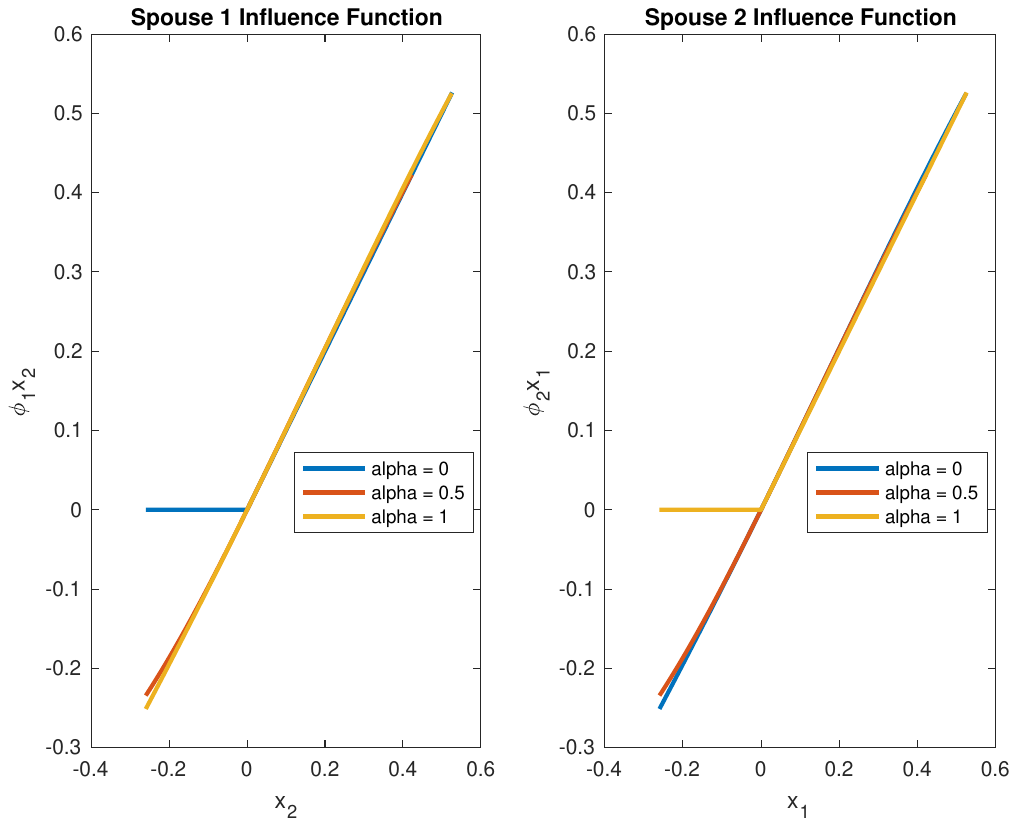}
    \caption{Sample Influence Functions for varying $\alpha$ values and fixed $\epsilon = 0.1$.}
    \label{fig:ex3}
\end{figure}

In figure \ref{fig:ex3}, we plot sample influence functions for different values of $\alpha$ with a fixed $\epsilon = 0.1$. When $\alpha = 0$, the couple is prioritizing the happiness of spouse 2, so spouse 1 adopts a conflict avoiding interaction style while spouse 2 is validating. On the other side of the spectrum, when $\alpha = 1$, the spouses are prioritizing the happiness of spouse 1, so spouse 1 adopts a validating interaction style while spouse 2 adopts a conflict avoiding. When $\alpha = 1/2$ and the couple is prioritizing their happiness equally, each spouse adopts a validating interaction style.

\section{Discussion}\label{discussion1}
Our research aims to recover the interaction styles identified by Gottman \cite{gottman_roles_1993}. In the process we showed that two of the interacton styles that couples adopt in successful marriages are optimal mathematical solutions when they have matched goals. These interaction styles were discovered from the research of Gottman \cite{gottman_marital_1989} observing how couples interact during an interview: Three styles are considered stable, and two are considered unstable, meaning these latter couples are more likely to get divorced. The three stable relationships are volatile, avoiding, and validating. Volatile couples are characterized by interactions in their negative range of emotions and few during their positive range of emotions. Validating couples interact in both their positive and negative emotions. Avoiders interact during their positive range of emotions, and few interact in their negative range. 

Our main result is that we are able to recover two of the stable interaction-styles through our methods and provide an explanation as how they are achieved.  

Our model does not produce Gottman's volatile interaction-style, which has been observed to be present in couples with the highest natural disposition. Although this style is not an optimal solution, it is nevertheless still capable of maintaining a stable marriage due to the high natural disposition of the couples. In other words, their marriage is stable despite their volatile interaction style.  However, adopting one of our optimal strategies would have led to a happier marriage for these couples.

\subsection{How robust is the model to terminal time $T$}

\begin{figure}[H]
    \centering
    \includegraphics[width=\textwidth, =\textwidth, trim= 2cm 2cm 2cm 2cm]{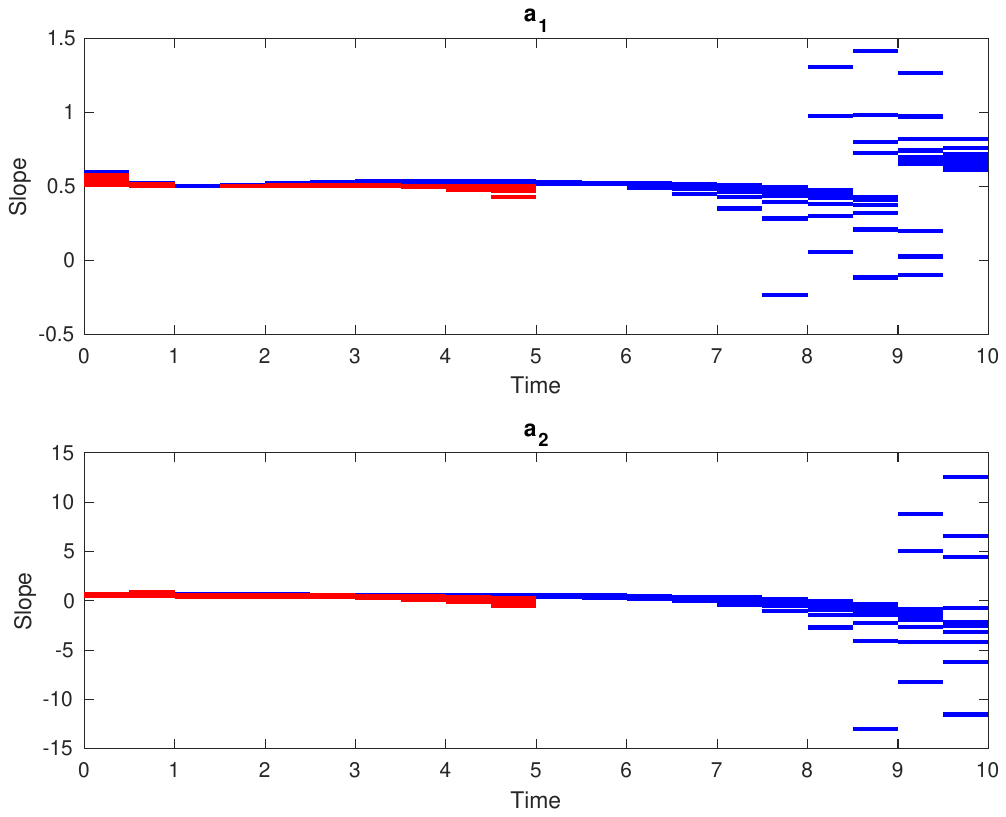}
    \caption{The slope averaged for time intervals of size 1 for terminal times of T = 5 and T = 10, for many different initial conditions}
    \label{fig:transInt}
\end{figure}

For this paper, one unit of time in our model corresponds to 5-minutes of real time (during the interview).  Gottman found that they could discover the couple's interaction styles during a 15- minute interview, and the same interaction styles are then used to predict the future of the marriage.  The concern was, could the result depend on the initial mood of the spouses when the interview began? And, if the interview was longer, would the researchers found a different influence function?. Next we test the robustness of our results with respect to $T$.  We could not make $T$ to be too large, because the numerical simulation does not produce reliable results for too large values of $T$; for that cases we can only verify numerical results with analytical results for small values of $\epsilon$. In Figure \ref{fig:transInt}, we plot the average slope of the influence function for different initial conditions for intervals of time $t = 1$ for terminal times of $T = 5$ (in blue) and $T = 10$ (in red). The terminal condition $\lambda(T)=0$ means that $u_i^*(T) = u_0$, regardless of the parameters. The simulation shows that the slope is fairly consistent up to $t = 3$ for a terminal time of $T= 5$ and $t = 6$ for a terminal time of $T = 10$. This means that regardless of whether the conversation is 15 minutes or 30 minutes, a similar interaction style is found, which leads to the same predicted outcome for the marriage. This is further evidence that Gottman's original research \cite{gottman_marital_1989} using short conversations to predict the success of a marriage is a valid assumption. 

Furthermore, Figure \ref{fig:transInt} shows a lot of variability in the slopes for times near the terminal time. The scatter is due to the fact that although the optimal control always collapses to $u_0$, different solutions approach it at different rate. 

\bibliography{citations}

\appendix
\section{Proof that $\lambda_i >0$}
\label{proof}
We prove this fact by exhausting all possible solutions of this problem. This can be broken into three cases:
\begin{enumerate}
    \item $u_i^* = u_{i,max}$
    \item $u_i^* = u_{i,max}, u_j^* = 0$
    \item $u_i^* = 0$,
\end{enumerate}
with $i,j = 1,2$ and $j\neq i$.

\subsection{Case 1: $u_i^* = u_{i,max}$}
To prove this case, we perform stability analysis on the adjoint system in equations (\ref{adj1}) and (\ref{adj2}). This system has an equilibrium point at
\begin{equation}\label{adj_EP}
    \left(\frac{(1-\alpha)u_2^* + \alpha r_2}{r_1r_2-u_1^*u_2^*}, \frac{(1-\alpha)r_1 + \alpha u_1^*}{r_1r_2-u_1^*u_2^*}\right).
\end{equation}

The adjoint system has characteristic polynomial $\mu^2 -(r_1+r_2)\mu + r_1r_2-ua$ = 0. The roots of the characteristic polynomial are 
\begin{equation} \label{adj_ew}
    \mu_{\pm} = \frac{r_1+r_2 \pm \sqrt{(r_1+r_2)^2- 4(r_1r_2-u_1^*u_2^*)}}{2}.
\end{equation}
The corresponding eigenvectors s are given by 
\begin{eqnarray} \label{adj_ev}
    v = \begin{bmatrix} \frac{u_2^*}{r_1-\mu_+} \\ 1 \end{bmatrix}, \\
    w = \begin{bmatrix} \frac{u_2^*}{r_1-\mu_-} \\ 1 \end{bmatrix}.
\end{eqnarray}
 
 There are three possible classifications of the equilibrium point. They are defined by what we call classification regions:
\begin{enumerate}
    \item  Where $u_1^*u_2^* > r_1r_2$, the equilibrium point is a saddle point since this is where $\mu_+ > 0$ and $\mu_- < 0$.
    \item Where $\frac{(r_1-r_2)^2}{-4}<u_1^*u_2^*<r_1r_2$, the equilibrium point is an unstable node because this is the region for which $\mu_{\pm}$ are positive.
    \item  Where $u_1^*u_2^* < \frac{(r_1-r_2)^2}{-4}$, the equilibrium point is an unstable spiral. This is where $\mu_{\pm}$ are complex.
\end{enumerate}

Because $u_{i,max} > 0$, we are only concerned with regions 1 and 2. 
\subsubsection{$u_{1,max}u_{2,max} > r_1r_2$ }
We have a saddle point. The numerator of the equilibrium point is always positive, but the sign of the denominator depends on the classification of the point. Here, the denominator is negative, so we can conclude that the point must be in the third quadrant. The eigenvalues given in equation (\ref{adj_ew}) have $\mu_+ >0$ and $\mu_-<0$. This tells us that eigenvector $v$ gives us the equation of the unstable manifold, and $w$ gives is the equation of the stable manifold. Therefore, $v$ has a negative slope, and $w$ has a positive slope. 

The slopes of the stable and unstable manifolds and the location of the equilibrium point in the third quadrant make it clear that any trajectory that will satisfy that transversality condition must originate in the first quadrant.

\subsubsection{$0 <u_{1,max}u_{2,max} < r_1r_2$} \label{unsNode}

These conditions tell us that the equilibrium point (\ref{adj_EP}) is an unstable node in the first quadrant. In order to satisfy the terminal conditions, $\lambda_1$ and $\lambda_2$ must begin in the first quadrant because beginning in any other quadrant would send the trajectories away from the first quadrant, and by extension also the origin.

\subsection{Case 2: $u_i^* = u_{i,max}, u_j^* = 0$}

Without loss of generality, assume $i = 1$ and $j = 2$. The equilibrium point for this system is now

\begin{equation}
    \left(\frac{\alpha}{r_1}, \frac{(1-\alpha)r_1 + \alpha u_1^*}{r_1r_2}\right).
\end{equation}
 The roots of the characteristic polynomial are $r_1$ and $r_2$ which are both positive and real. Thus, this equilibrium point is an unstable node. Just like in appendix \ref{unsNode}, $\lambda_1$ and $\lambda_2$ must be positive to satisfy the terminal condition. Therefore, $\lambda_1>0$ and $\lambda_2>0$.
\subsection{Case 3: $u_i^* = 0$}
The adjoint system is simplified to

\begin{subequations}
\begin{align}
        \dot{\lambda}_{1} &= -\alpha + \lambda_{1}r_1 ,\\
    \dot{\lambda}_{2} & = -(1-\alpha) + \lambda_{2}r2 
\end{align}
\end{subequations}
The equilibrium point is 
\begin{equation}\label{adj_EP2}
    \left(\frac{\alpha}{r_1}, \frac{1-\alpha}{r_2}\right).
\end{equation}
 This equilibrium point is in the first quadrant. $\dot{\lambda}_{1}$ is negative for any value of $\lambda_1 < \alpha/r_1$. Similarly, $\dot{\lambda}_{2}$ is decreasing for any value of $\lambda_2 < (1-\alpha)/r_2$. We come to the same conclusion. Any trajectory that satisfies the terminal condition must begin and remain in the first quadrant. Therefore, $\lambda_1>0$ and $\lambda_2>0$.
\section{Singular Solution}

Now, we must determine whether or not there is a singular solution. For $u_1^*$ to be singular, $\lambda_1x_2 \equiv 0$ for an interval of time. This means that $\lambda_1 = 0$ or $x_2 = 0$. First, we will consider the case for $\lambda_1 = 0$.

If $\lambda_1 = 0$, but $\dot{\lambda}_1 \neq 0$, then $u_1^*$ cannot be singular. Therefore, $\dot{\lambda}_1 = 0$ means
    \begin{eqnarray}
        0 = \dot{\lambda}_1 = -\alpha -\lambda_2u_2^*, \\
        \dot{\lambda}_2 = -(1-\alpha) + \lambda_2r_2,\\
        \text{ implies } \lambda_2(t) = \frac{1-\alpha}{r_2} -\frac{1-\alpha}{r_2}e^{r_2(t-T)},\\
        0 = -\alpha -u_2^*[ \frac{1-\alpha}{r_2}, -\frac{1-\alpha}{r_2}e^{r_2(t-T)}].
    \end{eqnarray}
    If $\alpha = 0$,
    \begin{eqnarray}
        0 = -u_2^*[ \frac{1}{r_2} -\frac{1}{r_2}e^{r_2(t-T)}],\\
        \frac{u_2^*}{r_2} = \frac{u_2^*}{r_2}e^{r_2(t-T)}, \\
        1 = e^{r_2(t-T)}, \\
        t = T.
    \end{eqnarray}
    This is not an interval.

    If $\alpha = 1$,
    \begin{eqnarray}
        0 = -1.\\
    \end{eqnarray}
    This is a contradiction. 
    If $\alpha \in (0,1)$,
    \begin{eqnarray}
        0 =  -\alpha -u_2^*[ \frac{1-\alpha}{r_2} -\frac{1-\alpha}{r_2}e^{r_2(t-T)}], \\
        u_2^* = \frac{\alpha}{-\frac{1-\alpha}{r_2} +\frac{1-\alpha}{r_2}e^{r_2(t-T)}}. \\
        \text{We need}
        0 <  \frac{\alpha}{-\frac{1-\alpha}{r_2} +\frac{1-\alpha}{r_2}e^{r_2(t-T)}} < M = u_{2,max} \\
       -\frac{1-\alpha}{r_2} +\frac{1-\alpha}{r_2}e^{r_2(t-T)} > 0\\
       e^{r_2(t-T)} > 1\\
       r_2(t-T) > 0\\
       t > T.\\
    \end{eqnarray}
This is a contradiction because $t \in [0,T]$. Therefore, $\lambda_1 \neq 0$ for an interval of time. 

For $u_2^*$ to be singular, $\lambda_2x_1\equiv 0$. The analysis for $\lambda_2$ is similar to the analysis for $\lambda_1$, so we come to the same conclusion that $\lambda_2 \neq 0$ for an interval of time. The remaining possible cases are

\begin{enumerate}
    \item $x_2 = 0, x_1 = 0$ and $u_1^*, u_2^*$ both singular, 
    \item $x_2 = 0$, $u_1^*$ singular, and $u_2^* = 0$ or $u_2^* = u_{2,max} = M_2$, and
    \item $x_1 = 0$, $u_2^*$ singular, and $u_1^* = 0$ or $u_1^* = u_{1,max} = M_1$.
\end{enumerate}

We will examine these cases one at a time. 

\begin{enumerate}
    \item Case: $x_2 = 0, x_1 = 0$ and $u_1^*, u_2^*$ both singular \\
    The state equations are reduced to
    \begin{subequations}
    \begin{align}
        0 = r_1\bar{x}_1 \\
        0 = r_2\bar{x}_2.
    \end{align}
    \end{subequations}
    The parameters $r_1$ and $r_2$ are never equal to zero. We are not interested in the cases where $\bar{x}_1$ or $\bar{x}_2$ are equal to zero. 

    \item Case: $x_2 = 0$, $u_1^*$ singular, and $u_2^* = 0$ or $u_2^* = u_{2,max} = M_2$
    This means:
    \begin{eqnarray}
        \dot{x}_1 = r_1(\bar{x}_1 - x_1) \\
        \dot{x}_2 = r_2\bar{x}_2 + u_2^*x_1.
    \end{eqnarray}
    Thus, we can solve for $x_1$.
    \begin{equation}
         x_1(t) = \bar{x}_1 + (x_1^0 - \bar{x}_1)e^{-r_1t}
    \end{equation}
   Plugging this into the equation for $x_2$,
   \begin{equation}
       \dot{x}_2 = r_2\bar{x}_2 + u_2^*[\bar{x}_1 + (x_1^0 - \bar{x}_1)e^{-r_1t}].
   \end{equation}
   This is not true for an interval if $\dot{x}_2 \neq 0$ at $x_2 = 0$. Therefore, $\dot{x}_2 = 0$ 

   \begin{eqnarray}
       0 =  r_2\bar{x}_2 + u_2[\bar{x}_1 + (x_1^0 - \bar{x}_1)e^{-r_1t}]  \label{x2eq0}\\ 
       u_2^* = \frac{-r_2\bar{x}_2}{\bar{x}_1 + (x_1^0 - \bar{x}_1)e^{-r_1t}} = M_2.
   \end{eqnarray}

For most parameters, equation (\ref{x2eq0}) only applies to specific parameters. However, if $x_1^0 = \bar{x}_1$ and
\begin{equation}
    -r_2\bar{x}_2 = M_2\bar{x}_1,
\end{equation}
$u_1^*$ will be a singular solution while $u_2^*$ is not.

\item Case: $x_1 = 0$, $u_2^*$ singular, and $u_1^* = 0$ or $u_1^* = u_{1,max} = M_1$\\
This is similar to the previous case. If $x_2^0 = \bar{x}_1$ and 
\begin{equation}
     -r_1\bar{x}_1 = M_1\bar{x}_2,
\end{equation}
$u_2^*$ will be singular while $u_1^*$ is not.
\end{enumerate}

\end{document}